\documentclass[12pt,centertags]{amsart}
\usepackage{amsmath,amstext,amsthm,a4,amssymb,amscd}
\usepackage[mathscr]{eucal}
\usepackage{mathrsfs}

\numberwithin{equation}{section}

\newtheorem{thm}{Theorem}[section]

\theoremstyle{definition}
\newtheorem{rem}[thm]{Remark}
\theoremstyle{definition}

\newcommand{\be}{\begin{eqnarray}}
\newcommand{\ee}{\end{eqnarray}}

\newcommand{\comment}[1]{}

\begin{document}
\title{The Mathematical Work of V. K. Patodi}

\author{Weiping Zhang}

\begin{abstract}
We give a brief survey on   aspects  of the local index theory as developed from   the mathematical works of V. K. Patodi.  It is dedicated to the 70th anniversary of Patodi. 
\end{abstract}

\maketitle

\setcounter{section}{-1}
\section{Introduction} \label{s0}

Vijay Kumar Patodi was born on March 12, 1945 in Guna, India. He passed away on December 21, 1976, at the age of 31. 

It is remarkable that even in such a short period of life, Patodi made quite a number of fundamental  contributions to mathematics.  These contributions have had deep and wide spread influences in geometry, topology, number theory, as well as in mathematical physics.  

Patodi gave an invited talk at the International Congress of Mathematicians in 1974. His Collected Papers \cite{P0}, edited by Atiyah and Narasimhan,  was published by World Scientific in 1996.

Patodi's work centers around the  Atiyah-Singer  index theory, a major subject in global analysis. He published 13 mathematical papers, all included in his Collected Papers \cite{P0}. Many of these papers have become classical literatures in mathematics.   They cover   almost all aspects of the index theory, from the classical period where the topological methods play a major role, to the modern era where more geometrical  methods, including those developed by Patodi himself, become more and more important.

In this article, we will give a survey  of Patodi's works on local index theorems, on $\eta$ invariant, and on analytic torsion, as well as their implications. Due to the limitation of the knowledge of the author, we only concentrate on subjects  we are most familiar. 


$
\
$

\noindent{\bf Acknowledgements.} {The author is partially supported by the National Natural Science Foundation
of China.}   He is grateful to Xianzhe Dai and Xiaonan Ma for very helpful suggestions.

\section{Local index theorem}\label{s1}

The celebrated index theorem proved by Atiyah and Singer \cite{AS} has been regarded as one of the most significant mathematical results of the last century. It connects many areas of mathematics. Roughly  speaking, it states that the analytically defined index of an elliptic differential operator on a closed  manifold can be computed by purely topologically defined quantities. Many famous results in geometry, topology and algebraic geometry are special examples of the Atiyah-Singer index theorem. We mention typically the Gauss-Bonnet-Chern theorem \cite{Ch45}, the Hirzebruch signature theorem \cite{H}, the Riemann-Roch-Hirzebruch  theorem for algebraic varieties  \cite{H} and the index theorem for Dirac operators on spin manifolds \cite{AS}. In order to have a thorough understanding of Patodi's first paper \cite{P1}, we take the Gauss-Bonnet-Chern theorem as a typical example. 

Let $M$ be a closed oriented manifold of   dimension $2n$. Let $g^{TM}$ be a Riemannian metric on the tangent vector bundle $TM$ of $M$. Let $\nabla^{TM}$ be the Levi-Civita connection on $TM$ associated to $g^{TM}$. Let $R^{TM}=(\nabla^{TM})^2$ be the curvature of $\nabla^{TM}$. This curvature locally can be viewed as a skew-adjoint matrix with elements consisting of 2-forms. For example, if $e_1,\,\cdots,\, e_{2n}$ is a (local) oriented orthonormal basis of $TM$ and $e^i$, $1\leq i\leq 2n$, is the dual basis of $T^*M$, the cotangent vector bundle of $M$, then each element $R^{TM}_{ij}$, with $1\leq i,\, j\leq 2n$, can be written as $R^{TM}_{ij}=\frac{1}{2}\sum_{s,\, t=1}^{2n}R^{TM}_{ijst}e^s\wedge e^t$. 

We define ${\rm Pf}(R^{TM})$, the Pfaffian of $R^{TM}$, by 
\begin{align}\label{1.1}
{\rm Pf}\left(R^{TM}\right)=\frac{1}{2^{ {n} }  {n} !}\sum_{i_1,\, \cdots,\, i_{2n}=1}^{2n}\epsilon_{i_1\cdots i_{2n}}R^{TM}_{i_1i_2}\wedge\cdots\wedge R^{TM}_{i_{2n-1}i_{2n}},
\end{align}
where $\epsilon_{i_1\cdots i_{2n}}$ is the permutation  number of $\{i_1,\,\cdots,\, i_{2n}\}$ with respect to $\{1,\,\cdots,\,2n\}$. 

Let $\chi(M)$ be the Euler characteristic number of $M$. 

The Gauss-Bonnet-Chern theorem \cite{Ch45} states that the following identity holds,
\begin{align}\label{1.2}
\chi(M)=\left(\frac{-1}{2\pi}\right)^n\int_M{\rm Pf}\left(R^{TM}\right) .
\end{align}
In particular, when $n=1$,   one gets the Gauss-Bonnet formula $\int_MKdv_M=2\pi\chi(M)$, where $K$ is the Gauss curvature of $g^{TM}$ and $dv_M$ is the volume form associated to $g^{TM}$. 

We first show how (\ref{1.2}) can be viewed as an index theorem in the sense of Atiyah and Singer. 

Let $\Lambda^*(T^*M)=\oplus_{i=0}^{2n}\Lambda^i(T^*M)$ be the exterior algebra bundle of $T^*M$. Let $\Omega^i(M)=\Gamma(\Lambda^i(T^*M) )$ be the space of smooth $i$-forms on $M$. Let $(\Omega^*(M),d)$ be the de Rham complex
\begin{align}\label{1.3}
\left(\Omega^*(M),d\right): 0\rightarrow\Omega^0(M)\stackrel{d}{\rightarrow}\Omega^1(M)\stackrel{d}{\rightarrow}\cdots\stackrel{d}{\rightarrow}\Omega^{2n}(M)\rightarrow 0.  
\end{align}

The metric $g^{TM}$ induces natural metrics  on each $\Lambda^i(T^*M)$ and thus on $\Lambda^*(T^*M)$ through orthogonal direct sum. They induce a natural inner product $\left\langle \cdot,\cdot\right\rangle$ on $\Omega^*(M)$. 

Let $d^*$ be the formal adjoint of $d$ with respect to the inner product on $\Omega^*(M)$. That is, for any $\alpha,\, \beta\in \Omega^*(M)$, one has $\langle d\alpha,\beta\rangle=\langle \alpha,d^*\beta\rangle$. The operator
\begin{align}\label{1.4}
D=d+d^*:\Omega^*(M)\longrightarrow \Omega^*(M)
\end{align}
is called the de Rham-Hodge operator associated to $g^{TM}$. It is a formally self-adjoint elliptic differential operator.  Let 
\begin{align}\label{1.5}
D_+=d+d^*:\Omega^{\rm even}(M)=\bigoplus_{i=0}^n \Omega^{2i}(M)\longrightarrow \Omega^{\rm odd}(M)=\bigoplus_{i=0}^{n-1}\Omega^{2i+1}(M),
\end{align}
$$D_-=d+d^*:\Omega^{\rm odd}(M) \longrightarrow \Omega^{\rm even}(M) $$
be the natural restrictions of $D$. Then $D_\pm$ are  elliptic differential operators and $D_-$ is the formal adjoint of $D_+$. Moreover, 
by the Hodge decomposition theorem, one finds that the index of $D_+$ equals to $\chi(M)$: 
\begin{align}\label{1.6}
 \chi(M)={\rm ind}\left(D_+\right):=\dim\left(\ker D_+\right)-\dim\left(\ker D_-\right).
\end{align}

By (\ref{1.6}), one can reformulate  the Gauss-Bonnet-Chern formula (\ref{1.2}) as an index formula
\begin{align}\label{1.7}
 {\rm ind}\left(D_+\right) =\left(\frac{-1}{2\pi}\right)^n\int_M{\rm Pf}\left(R^{TM}\right) .
\end{align}

This is a typical form of the Atiyah-Singer index theorem. For a different elliptic differential operator $D_+$, the right hand side will have different topological meaning. 

We now explain a simple and beautiful idea due to Mckean and Singer \cite{MS}. 

Let $E$, $F$ be two Hermitian vector bundles over $M$. Let $D_+:\Gamma(E)\rightarrow \Gamma(F)$ be a first order elliptic differential operator. Let $D_-:\Gamma(F)\rightarrow \Gamma(E)$ be the formal adjoint of $D_+$. Then both $D_-D_+:\Gamma(E)\rightarrow \Gamma(E)$ and $D_+D_-:\Gamma(F)\rightarrow \Gamma(F)$ are nonnegative formally self-adjoint elliptic operators. In particular, for any $t>0$, the heat operators $\exp(-tD_-D_+)$ and $\exp(-tD_+D_-)$ are well-defined. 

The following simple formula due to Mckean-Singer is the starting point of the whole local index theory: 
\begin{multline}\label{1.8}
 {\rm ind}\left(D_+\right) ={\rm Tr}\left[\exp\left(-tD_-D_+\right)\right]- {\rm Tr}\left[\exp\left(-tD_+D_-\right)\right]\\
=\sum_{\lambda\in{\rm Spec}(D_-D_+)}\exp(-t\lambda)-\sum_{\mu\in{\rm Spec}(D_+D_-)}\exp(-t\mu) ,
\end{multline}
where ``Spec$(\cdot)$'' denotes the set of spectrum (or eigenvalues) of the corresponding elliptic operator. 

The proof of  (\ref{1.8}) is very simple: if $\lambda\in {\rm Spec}(D_-D_+)$ and $\lambda\neq 0$, then by $D_-D_+s=\lambda s$ with $s\neq 0$, one gets $D_+D_-(D_+s)=\lambda(D_+s)$ with $D_+s\neq 0$. Thus $\lambda\in {\rm Spec}(D_+D_-)$. Similarly, any  nonzero $\mu\in {\rm Spec}(D_+D_-)$ belongs also to $ {\rm Spec}(D_-D_+)$. On the other hand,   one verifies that  $\ker(D_\mp D_\pm)=\ker(D_\pm)$. In summary, one finds
\begin{multline}\label{1.9}
 \sum_{\lambda\in{\rm Spec}(D_-D_+)}\exp(-t\lambda)-\sum_{\mu\in{\rm Spec}(D_+D_-)}\exp(-t\mu)=\dim(\ker D_+)-\dim(\ker D_-)\\ ={\rm ind}\left(D_+\right) ,
\end{multline}
from which (\ref{1.8}) follows. 

Formula (\ref{1.8}) suggests new interpretations of the index.  In fact, as both $\exp(-tD_-D_+):\Gamma(E)\rightarrow \Gamma(E)$ and $\exp(-tD_+D_-):\Gamma(F)\rightarrow \Gamma(F)$ are easily seen to be compact operators, they admit smooth  kernels (called heat kernels) $P_t(\cdot,\cdot)$, $Q_t(\cdot,\cdot)$  respectively such that for any $x,\, y\in M$, $P_t(x,y)\in {\rm Hom}(E_y,E_x)$, $Q_t(x,y)\in {\rm Hom}(F_y,F_x)$, and that for any $u\in\Gamma(E)$, $v\in\Gamma(F)$, one has
\begin{align}\label{1.10}
\left(\exp\left(-tD_-D_+\right) u\right)(x)=\int_MP_t(x,y)u(y)dv_M(y),
\end{align}
$$ \left(\exp\left(-tD_+D_-\right) v\right)(x)=\int_MQ_t(x,y)v(y)dv_M(y).$$

From (\ref{1.8}) and (\ref{1.10}), one gets
\begin{align}\label{1.11}
{\rm ind}\left(D_+\right)=  \int_M\left({\rm tr}\left[P_t(x,x)\right]  -  {\rm tr}\left[Q_t(x,x)\right]\right) dv_M(x).
\end{align}

The nice point of (\ref{1.11}) is that since the left hand side does not depend on $t>0$, one can deform $t$ in the right hand side to try to get more information about the left hand side. Indeed, the asymptotic behavior as $t\rightarrow 0$ of $P_t(\cdot,\cdot)$ and $Q_t(\cdot,\cdot)$ has been studied in classical literature (cf. \cite{MP}). In particular, there exist smooth functions $a_0,\, \cdots,\, a_n$ (resp. $b_0,\,\cdots,\, b_n$) on $M$ such that when $t\rightarrow 0^+$, one has for any $x\in M$,
\begin{align}\label{1.12}
{\rm tr}\left[P_t(x,x)\right]  =\frac{1}{(4\pi t)^n}\left(a_0(x)+a_1(x)\,t+\cdots +a_n(x)\,t^n+O\left(t^{n+1}\right)\right),
\end{align}
$$ {\rm tr}\left[Q_t(x,x)\right]  =\frac{1}{(4\pi t)^n}\left(b_0(x)+b_1(x)\,t+\cdots +b_n(x)\,t^n+O\left(t^{n+1}\right)\right).$$

Combining (\ref{1.11}) with (\ref{1.12}), one gets
\begin{align}\label{1.13}
{\rm ind}\left(D_+\right)= \sum_{i=0}^n\frac{t^i}{(4\pi t)^n} \int_M\left( a_i(x)-b_i(x)\right) dv_M(x)+O(t).
\end{align}

Since $t>0$ can vary, while the left hand side of (\ref{1.13}) does not depend on $t$, from (\ref{1.13}) one deduces that if $0\leq i\leq n-1$, then
\begin{align}\label{1.14}
\int_M\left( a_i(x)-b_i(x)\right) dv_M(x)=0,
\end{align}
while if $i=n$, one has
\begin{align}\label{1.15}
\frac{1}{(4\pi)^n}\int_M\left( a_n(x)-b_n(x)\right) dv_M(x)={\rm ind}\left(D_+\right).
\end{align}
Combining with (\ref{1.7}), one gets, in the case of $D_+=d+d^*$, that
\begin{align}\label{1.15a}
\frac{1}{(4\pi)^n}\int_M\left( a_n(x)-b_n(x)\right) dv_M(x)=\left(\frac{-1}{2\pi}\right)^n\int_M{\rm Pf}\left(R^{TM}\right).
\end{align}

Now at least for the de Rham-Hodge operator $D_\pm=d+d^*$, by \cite{MP} one knows that  theoretically,  the functions $a_i$'s and  $b_i$'s can be expressed by local (pointwise) geometric quantities associated directly to the Riemannian metric $g^{TM}$. It leads Mckean and Singer to pose in \cite{MS} the conjecture that for the classical geometric operators such as the de Rham-Hodge operator, one can refine (\ref{1.14}) and (\ref{1.15}) so that for any $0\leq i\leq n-1$, one has\footnote{That $a_0=b_0=1$ is a classical result of Hermann Weyl.}
\begin{align}\label{1.16}
 a_i -b_i =0
\end{align}
over $M$, while for $i=n$, $\frac{1}{(4\pi)^n} \left( a_n(x)-b_n(x)\right)dv_M(x) $, which we  call the {\it index density},  would provide  the geometric/topological information desired  by the Atiyah-Singer index theorem. In particular, for the de Rham-Hodge operator $D_+=d+d^*:\Omega^{\rm even}(M)\rightarrow \Omega^{\rm odd}(M)$, one may expect 
\begin{align}\label{1.17}
\frac{1}{(4\pi)^n} \left( a_n(x)-b_n(x)\right) dv_M(x) = \left(\frac{-1}{2\pi}\right)^n {\rm Pf}\left(R^{TM}\right) .
\end{align}

Mckean and Singer were able to prove (\ref{1.17}) for the case of $n=1$, while in this case (\ref{1.16}) is a trivial consequence of the Hermann Weyl result. They went on to raise (\ref{1.16}) and (\ref{1.17}) for the de Rham-Hodge operator as a conjecture of ``fantastic cancellation''.

We should point out that while theoretically the functions $a_i$'s and $b_i$'s can be constructed from $g^{TM}$, the process is very complicated  and involves taking higher order derivatives of $g^{TM}$. This is why Mckean and Singer called their conjecture, if it would hold, {\it fantastic}. 

We can now describe the contribution obtained by Patodi in his very first paper \cite{P1}. If one should  describe it in one sentence, then one simply says that in \cite{P1}, Patodi proved that (\ref{1.16}) and (\ref{1.17}) hold for the de Rham-Hodge operator. That is, the ``fantastic cancellation'' expected by Mckean and Singer does hold!

It is easy to see that by integrating (\ref{1.16}) and (\ref{1.17}) over $M$, we  get (\ref{1.2}). Thus, (\ref{1.16}) and (\ref{1.17}) together refine (\ref{1.2}) significantly. It is a local pointwise formula. For this reason we call it a local index theorem, and in this case we can say that Patodi established the local Gauss-Bonnet-Chern theorem. 

The proof of Patodi in \cite{P1} for (\ref{1.16}) and (\ref{1.17}) is a true {\it tour de force}. From the work of Minakshisundaram and   Pleijel \cite{MP}, when $t>0$ is very small, one can approximate  $P_t(x,y)$ and $Q_t(x,y)$ by the so called paramatrices which are determined by a series of inductive equations involving the Laplacians $D_-D_+$ and $D_+D_-$, as well as other   local geometric data. It takes remarkable insight and highly nontrivial computations to show that they can all be put into an ordered manner so that the final contribution indeed gives (\ref{1.16}) and (\ref{1.17}). 

While the paper \cite{P1} already gives a breakthrough contribution by showing that the local index theorem does hold for the de Rham-Hodge operator, in his second paper \cite{P2}, Patodi moved on to show that the local index theorem also holds for the Dolbeault operator on K\"ahler manifolds, whose global version implies the Riemann-Roch-Hirzebruch theorem for algebraic manifolds \cite{H}. He also mentioned in the end of \cite{P2} that he was able to prove the local index theorem for the signature operator whose global version is the Hirzebruch signature theorem \cite{H}. 

The computations    in \cite{P2} are   more complicated than those given in \cite{P1}, and it is truly remarkable that Patodi was able to accomplish all these while staying   in India. 

Mathematically, among the four theorems mentioned above, now only the local index theorem for Dirac operator was missing, and it is the subject of the joint paper \cite{P4} with Atiyah and Bott. 

Before describing the content of \cite{P4}, we mention that almost during  the same time as Patodi established local index theorems for the de Rham-Hodge and Dolbeault operators by direct tour de force methods, Gilkey, on the other hand, developed an indirect way to prove   local index theorems in \cite{Gilkey73a}, \cite{Gilkey73b}, by using Weyl's invariant theory. 

In \cite{P4}, together with Atiyah and Bott, Patodi examined the local index theorem systematically by combining his direct approach  with Gilkey's method. The authors of \cite{P4} showed that the local index theorem holds for the so called twisted signature operator, which allows to twisting  a vector bundle to the signature operator,  as well as the Dirac   and twisted Dirac operators. Moreover, from the twisted signature theorem, by combining with the Bott periodicity theorem, they arrived at a new proof of the (global) Atiyah-Singer index theorem for general elliptic differential operators.

There are two points one should mention about \cite{P4}. The first is that the method and result in \cite{P4} enable one to generalize the  index theorem for geometric operators such as the signature operator and the Dirac operator to the case of manifolds with boundary. This is the content of \cite{P6} and \cite{P7} and will be described in the next section. The other one is that while the local Gauss-Bonnet-Chern and the local Riemann-Roch-Hirzebruch theorems can be proved directly in \cite{P1} and \cite{P2}, the proof of the local index theorem for Dirac operators in \cite{P4} is not direct. The direct proof of the local index theorem for Dirac operators, by using the method of Patodi, was later given by   Yu \cite{Yu87}.

In \cite{P5} and in the joint paper with Donnelly \cite{P11}, Patodi generalized his method to the equivariant situation and proved the Lefschetz fixed point formula   of Atiyah-Bott-Segal-Singer (\cite{AB}, \cite{ASe}, \cite{ASI}) by the heat equation method.  The paper \cite{P5} is an announcement for the holomorphic case, while the paper \cite{P11} covers the case of the equivariant signature theorem with respect to an isometry on a Riemannian manifold. 

We would like to mention that locally, all the geometric operators mentioned above: the de Rham-Hodge operator, the signature operator and the Dolbeault operator on K\"ahler manifolds can be expressed as a twisted Dirac operator.  So now by local index theorem one usually means the local index theorem for (twisted) Dirac operators.  

In early 1980's, several physics inspired formal proofs of the index theorem for Dirac operators   appeared. Among them, 
Atiyah described in \cite{A85} an idea due to the physicist  Witten, who suggested that by formally applying the Duistermaat-Heckman  localization formula \cite{DH82} in equivariant cohomology to the loop space of a closed spin manifold, one can obtain the index theorem for Dirac operators tautologically, though non-rigorously. The ideas of Witten and other physicists (cf. \cite{Al} and \cite{FW84}) stimulated a number of new proofs of the local index theorem for Dirac operators, notably the probabilitistic proof due to Bismut \cite{B84}, the two proofs of Getzler (\cite{G83}, \cite{G86})  inspired by supersymmetry and the group-theoretic proof given by Berline-Vergne \cite{BV}.  It turns out that the index density is closely related to the heat kernel of the harmonic oscillator. 

The proofs of Bismut \cite{B84} and Berline-Vergne \cite{BV} also cover the equivariant case. Another heat kernel proof of the equivariant index theorem for Dirac operators, which is closely related to the method developed by Patodi (which is  extended to   spin manifolds   in \cite{Yu87}), is given in \cite{LYZ}. 

It is amazing that Patodi's direct computation now has connections with so many areas: from physics to probability, and to group representations. 

In ICM-1986, Bismut gave an invited talk entitled ``Index theorem and the heat equation'' \cite{B86a}.  Besides providing an overview of the history of the heat equation approach to the index theorem, which was started from the pioneering works of Mckean-Singer and Patodi, the talk also describes the most recent advances due to the speaker himself: a vast  generalization on the heat equation proofs of the   index theorem   for a family of Dirac operators \cite{B86}, by making use of Quillen's concept of superconnection \cite{Q}. It marks the beginning  of a new era of the local index theory. 

We refer the interested reader  to the book of Berline-Getzler-Vergne \cite{BGV} for a systematic  treatment of the local index theorems for geometric operators, as well as the Bismut local families index theorem for Dirac operators \cite{B86}.

\section{$\eta$-invariant and the index theorem on manifolds with boundary} \label{s2} 
It is always important in analysis to study differential equations with boundary conditions. Thus even at the early stage of the development of the index theory, the problem of index theorems for elliptic operators on manifolds with boundary was studied, cf. \cite{AB64}. However, among the four geometric operators we have mentioned before: the de Rham-Hodge operator, the signature operator, the Dolbeault operator  on K\"ahler manifolds and the Dirac operator on spin manifolds, only the de Rham-Hodge operator admits the elliptic boundary condition in the classical context    studied in \cite{AB64}. These boundary conditions are ``local'' in the sense that if $P$ is such a boundary condition which acts on smooth sections on the boundary, then for any   smooth section $u$ on the boundary, if $u$ vanishes on an open subset of the boundary,   one requires that $Pu$ also vanishes on this open subset. It is known that the result in \cite{AB64} does not apply to   signature and/or Dirac operators. 

The problem of the index theorem for signature and Dirac operators on manifolds with boundary was solved in a series of joint papers \cite{P7} due to Atiyah, Patodi and Singer, with an earlier  announcement presented  in \cite{P6}.  
In this section, we will discuss this solution and its ramifications. 

This section consists of  two subsections. In Subsection \ref{s2.1}, we present the index theorem. In Subsection \ref{s2.2}, we will discuss in more details the $\eta$-invariant which occurs as the contribution from  the boundary in the index theorem. 

\subsection{The Atiyah-Patodi-Singer index theorem}\label{s2.1}

It is discovered in \cite{P7} that for a Dirac operator on a spin manifold  with boundary to be elliptic, one needs to impose a {\it global} boundary condition, instead of the usual {\it local} ones. 

To be more precise, let $M$ be an oriented spin manifold with smooth boundary $\partial M$.  We assume $M$ is of even dimension $2n$, then $\partial M$ is of dimension $2n-1$, and carries an induced orientation and also an induced spin structure. Let $g^{TM}$ be a Riemannian metric on $TM$. We assume that $g^{TM}$ is of {\it product structure} near the boundary $\partial M$. That is, there exists a neighborhood $\partial M\times [0,a)\subset M$ of $\partial M$, with $a>0$ sufficiently small, such that 
\begin{align}\label{2.1}
\left.g^{TM}\right|_{\partial M\times [0,a)}=\pi^*\left(\left.g^{TM}\right|_{\partial M}\right)\oplus dt^2,
\end{align}
where $t\in[0,a)$ is the parameter and $\pi:\partial M\times [0,a)\rightarrow \partial M=\partial M\times\{0\}$ denotes the canonical projection. 

Let $S(TM)=S_+(TM)\oplus S_-(TM)$ be the Hermitian vector bundle of spinors associated to $(TM,g^{TM})$.  Let $(E,g^E)$ be a Hermitian vector bundle over $M$ carrying a Hermitian connection $\nabla^E$. We assume that $g^E$ and $\nabla^E$ are of {\it product structure} on $\partial M\times [0,a)$. That is, 
\begin{align}\label{2.2}
\left.g^E\right|_{\partial M\times [0,a)}=\pi^*\left(\left.g^E\right|_{\partial M}\right),\ \ \  \left.\nabla^E\right|_{\partial M\times [0,a)}=\pi^*\left(\left.\nabla^E\right|_{\partial M}\right).
\end{align}

Let $D^E:\Gamma(S(TM)\otimes E)\rightarrow \Gamma(S(TM)\otimes E)$ denote the twisted (by $E$) Dirac operator defined by the above geometric data. Let $D^E_\pm:\Gamma(S_\pm(TM)\otimes E)\rightarrow \Gamma(S_\mp(TM)\otimes E)$ denote the obvious restrictions of $D^E$. 

In view of (\ref{2.1}) and (\ref{2.2}), one finds that on $\partial M\times [0,a)$, one has
\begin{align}\label{2.3}
D^E=c\left(\frac{\partial}{\partial t}\right)\left(\frac{\partial}{\partial t}+\pi^*D^E_{\partial M}\right),
\end{align}
where $c(\cdot)$ is the notation for the Clifford action and $D^E_{\partial M}:\Gamma((S(TM)\otimes E)|_{\partial M})\rightarrow \Gamma((S(TM)\otimes E)|_{\partial M})$ is the induced Dirac operator on $\partial M$. Let $D^E_{\partial M,\pm}:\Gamma((S_\pm(TM)\otimes E)|_{\partial M})\rightarrow \Gamma((S_\pm(TM)\otimes E)|_{\partial M})$ denote the obvious restrictions. Then both $D^E_{\partial M,\pm}$ are elliptic and formally self-adjoint. We also call them the induced Dirac operators on the boundary $\partial M$. 

Let $L^2((S_\pm(TM)\otimes E)|_{\partial M})$ be the $L^2$-completions of $\Gamma((S_\pm(TM)\otimes E)|_{\partial M})$. Let $L^2_{\geq 0}((S_\pm(TM)\otimes E)|_{\partial M})\subset L^2((S_\pm(TM)\otimes E)|_{\partial M})$ and  $L^2_{>0}((S_\pm(TM)\otimes E)|_{\partial M})\subset L^2((S_\pm(TM)\otimes E)|_{\partial M})$ be defined by 
\begin{align}\label{2.4}
 L^2_{\geq 0}((S_\pm(TM)\otimes E)|_{\partial M})=\bigoplus_{\lambda\in{\rm Spec}\left(D^E_{\partial M,\pm}\right),\, \lambda\geq 0}E_\lambda,
\end{align}
$$ L^2_{>0}((S_\pm(TM)\otimes E)|_{\partial M})=\bigoplus_{\lambda\in{\rm Spec}\left(D^E_{\partial M,\pm}\right),\, \lambda> 0}E_\lambda,$$
where $E_\lambda$ denotes the eigenspace of $\lambda$. Let $P^E_{\pm,\geq 0}$, $P^E_{\pm,>0}$ denote the orthogonal projections from 
$ L^2 ((S_\pm(TM)\otimes E)|_{\partial M})$ to $ L^2_{\geq 0}((S_\pm(TM)\otimes E)|_{\partial M})$, $ L^2_{> 0}((S_\pm(TM)\otimes E)|_{\partial M})$ respectively. 

The first result in \cite{P7} shows that the boundary problems $(D_\pm^E,P^E_{\pm,\geq 0})$  and $(D_\pm^E,P^E_{\pm,> 0})$ are elliptic. Moreover, $(D_-^E,P^E_{-,> 0})$  is the adjoint of  $(D_+^E,P^E_{+,\geq  0})$.  Thus one can define the index of $(D_+^E,P^E_{+,\geq  0})$ by
\begin{align}\label{2.5}
 {\rm ind}\left(D_+^E,P^E_{+,\geq  0}\right)=\dim \left(\ker \left(D_+^E,P^E_{+,\geq  0}\right)\right)-\dim \left(\ker \left(D_-^E,P^E_{-,>  0}\right)\right),
\end{align}
where both
\begin{align}\label{2.6}
 \ker \left(D_+^E,P^E_{+,\geq  0}\right) =\left\{ u\in \Gamma\left(S_+(TM)\otimes E\right):D^E_+u=0,\ P^E_{+,\geq  0}\left(u|_{\partial M}\right)=0\right\}
\end{align}
and 
\begin{align}\label{2.7}
 \ker \left(D_-^E,P^E_{-,>  0}\right) =\left\{ u\in \Gamma\left(S_-(TM)\otimes E\right):D^E_-u=0,\ P^E_{-,>  0}\left(u|_{\partial M}\right)=0\right\}
\end{align}
are of finite dimension. 

One of the immediate observations is that here $P^E_{\pm,\geq 0}$ and $P^E_{\pm,>0}$ no longer verifies the requirements  of the local boundary condition. Thus they are {\it global} boundary conditions. The price one pays  here is that these boundary conditions no longer contribute topologically  invariant indices. In particular, $ {\rm ind}(D_+^E,P^E_{+,\geq  0} )$ defined in (\ref{2.5}) now depends on the induced  geometric data on $\partial M$.  

The first main result in \cite{P7} provides an explicit formula for  $ {\rm ind}(D_+^E,P^E_{+,\geq  0} )$, generalizing the Atiyah-Singer index theorem for $ {\rm ind}(D_+^E)$ in the case where $M$ is closed. The amazing thing here is that besides the usual geometric term which can be seen from the local index computation, there appears in this index formula a new term contributed from the boundary $\partial M$. Thus before stating the index formula, let us first describe this extra term: the $\eta$-invariant associated to the induced Dirac operator  $D^E_{\partial M,+}$ on $\partial M$. 

Recall that $D^E_{\partial M,+}:\Gamma((S_+(TM)\otimes E)|_{\partial M})\rightarrow \Gamma((S_+(TM)\otimes E)|_{\partial M})$ is elliptic and formally self-adjoint.  By standard elliptic theory, one sees that for any complex number $s$ with ${\rm Re}(s)>>0$, the following $\eta$-function of $D^E_{\partial M,+}$ is well-defined,
\begin{align}\label{2.8}
 \eta\left(D^E_{\partial M,+},s\right)=\sum_{\lambda\in {\rm Spec}\left(D^E_{\partial M,+}\right),\, \lambda\neq 0}\frac{{\rm sgn}(\lambda)}{|\lambda|^s},
\end{align}
where ${\rm sgn}(\lambda)=1$ if $\lambda>0$, while ${\rm sgn}(\lambda)=-1$ if $\lambda<0$.  It is a spectral function depending only on the restriction of the geometric data $(g^{TM},g^E,\nabla^E)$ on the boundary. 

It is shown in \cite{P7} that $\eta(D^E_{\partial M,+},s)$ can be extended to a meromorphic function of $s$ over ${\bf C}$, which is holomorphic at $s=0$. Naturally, one calls $\eta(D^E_{\partial M,+},0)$ the {\it $\eta$-invariant} of $D^E_{\partial M,+}$ and denotes it by $\eta(D^E_{\partial M,+})$. If one also counts the zero eigenvalue of $D^E_{\partial M,+}$, one   defines as in \cite{P7}
\begin{align}\label{2.9}
 \overline{\eta}\left(D^E_{\partial M,+}\right)= \frac{\dim\left(\ker D^E_{\partial M,+}\right)+  \eta\left(D^E_{\partial M,+} \right)}{2}
\end{align}
and calls it the {\it reduced} $\eta$-invariant of $D^E_{\partial M,+}$. 

One can now state the index theorem for $ (D_+^E,P^E_{+,\geq  0} )$ as follows. It is stated in \cite[I\,(4.3)]{P7}.
\begin{thm}\label{t2.1} The following identity holds,
\begin{align}\label{2.10}
{\rm ind}\left(D_+^E,P^E_{+,\geq  0} \right)=\int_M\widehat{A}\left(TM,\nabla^{TM}\right){\rm ch}\left(E,\nabla^E\right)- \overline{\eta}\left(D^E_{\partial M,+}\right) ,
\end{align}
where $\widehat{A} (TM,\nabla^{TM} )$ and ${\rm ch} (E,\nabla^E )$ are the Hirzebruch $\widehat{A}$-characteristic form associated to $\nabla^{TM}$ and the Chern character form associated to $\nabla^E$ defined respectively by
\begin{align}\label{2.11}
\widehat{A}\left(TM,\nabla^{TM}\right)={\rm det}^{ {1}/{2}}\left(\frac{R^{TM}/4\pi}{\sinh\left(R^{TM}/4\pi\right)}\right)
\end{align}
  and 
\begin{align}\label{2.12}
{\rm ch}\left(E,\nabla^E\right)={\rm tr}\left[\exp\left(\frac{\sqrt{-1}}{2\pi}\left(\nabla^E\right)^2\right)\right].
\end{align}
\end{thm}

\begin{rem}\label{t2.2}  If $M$ is closed, then $\partial M=\emptyset$ and (\ref{2.10}) becomes ${\rm ind}(D_+^E)=\langle \widehat{A}(TM){\rm ch}(E),[M]\rangle$, which is the original Atiyah-Singer index theorem for $D_+^E$ (\cite{AS}). 
\end{rem}
\begin{rem}\label{t2.3a}  If the metrics $g^{TM}$, $g^E$ and $\nabla^E$ are not of product structure near the boundary $\partial M$, then one can first deform them to ones of the product structure (without changing their restrictions on the boundary) and   apply Theorem \ref{t2.1} to get a geometric formula involving the new geometric data, and then obtain a formula involving only  the original geometric data by using the Chern-Weil theory.  For a systematic treatment of the index theorem for geometric operators without  the ``product structure'' near boundary assumption, see \cite{Gilkey93}. 
\end{rem}
\begin{rem}\label{t2.3} In its most general form, the main result of \cite[I]{P7}, stated in \cite[I, Theorem (3.10)]{P7},  holds for any first order   elliptic differential  operator on manifolds with boundary. However,   only for the geometric operators one gets the explicit local expression for the interior contribution. On the other hand, as we have mentioned, all geometric operators of interests can locally be expressed as a kind of (twisted) Dirac operators. This is why in the above we state   explicitly the index theorem for Dirac operators. 
\end{rem}

One of the motivations of \cite{P7} is to solve a conjecture of Hirzebruch \cite{H73} concerning the computation of the signature of Hilbert modular varieties  having cusp singularities. We mention here that by making use of the version of Theorem \ref{t2.1} for signature operators, this Hirzebruch conjecture was later solved  by Atiyah, Donnelly and Singer (all have collaborated with Patodi) in \cite{ADS}. An independent proof was given by M\"uller \cite{Mu84}.  These two works establish a solid place for Theorem \ref{t2.1} in number theory. We refer to \cite{A87} for an overview on this aspect.  

Now let us examine Theorem \ref{t2.1} from purely index theoretic point of view. 

One observes first that by (\ref{2.10}), $\overline{\eta} (D^E_{\partial M,+} )$ depends on $(g^{TM},g^E,\nabla^E)|_{\partial M}$. Thus it is not a topological invariant. Moreover, it does not admit a local expression, i.e., it can not be expressed as an integration of  local geometric terms. This later fact is in fact pointed out in \cite{P7}, by indicating that the $\eta$-invariant is not multiplicative under finite coverings.  

The second observation from (\ref{2.10}) is that when ${\rm mod}\ {\bf Z}$, $\overline{\eta} (D^E_{\partial M,+} )$ depends smoothly on $(g^{TM},g^E,\nabla^E)|_{\partial M}$. Moreover, by combining with the Chern-Weil theory, one may calculate explicitly the variation of this $\eta$-invariant with respect to   smooth deformations of $(g^{TM},g^E,\nabla^E)|_{\partial M}$. 

The third observation is that if $(g^{TM},g^E,\nabla^E)|_{\partial M}$ changes, then $\dim (\ker D^E_{\partial M,+})$ might jump, thus by (\ref{2.10}),  ${\rm ind} (D_+^E,P^E_{+,\geq  0}  )$ depends in general on $(g^{TM},g^E,\nabla^E)|_{\partial M}$. Indeed this can also be seen directly from (\ref{2.5})-(\ref{2.7}), as when $\dim (\ker D^E_{\partial M,+})$ jumps,    $ P^E_{+,\geq  0}   $  also jumps and certainly we get a different elliptic boundary value problem which may give  a possibly different index. 

The problem of characterizing  the variation of the index in Theorem \ref{t2.1}  with respect to different elliptic boundary conditions   was solved in \cite[III]{P7} by introducing the concept of {\it spectral flow}. 

Theorem \ref{t2.1}, as well as the associated concepts of $\eta$-invariant and spectral flow, have later on played very important roles in many aspects of geometry, topology, number theory (as mentioned above) and mathematical physics. We will describe  some of these ramifications in the next subsection. Here we just quote what Atiyah wrote in his Forward to  \cite{P0} concerning  \cite{P7}: ``It was a great collaboration, exploiting the different talents of the participants, and I am glad it came to a successful conclusion."  On another place,   Atiyah wrote in the ``Commentary on Papers on Index Theory'' in his Collected Works \cite{A88} about \cite{P7} that ``In many ways the papers on spectral asymmetry were perhaps the most satisfying ones I was involved with.'' 

We conclude this subsection by mentioning that for Theorem \ref{t2.1} itself, there are now several books and different proofs devoting  to it. For the books, we mention two: one is due to Boo{\ss}-Bavnbek and Wojciechowski \cite{BoW93} which gives a fairly complete treatment following the lines of \cite{P7}; the other one is due to Melrose \cite{Me} which gives a different treatment  on the analytic aspects of the problem. 
For other   proofs of Theorem \ref{t2.1}, we   refer to \cite[I]{BC90} and \cite{DZ00}.

\subsection{$\eta$-invariants and its applications}\label{s2.2}
We first point out that the $\eta$-invariant, or rather the reduced $\eta$-invariant defined in (\ref{2.9}), for (twisted) Dirac operators can in fact be defined for any  formally self-adjoint elliptic differential operator on a   closed manifold.  For odd dimensional manifold case, this has been studied systematically in \cite[III, Sections 2-4]{P7}. For brevity, we will concentrate on   Dirac operators. 

For any closed odd dimensional  oriented Riemannian spin manifold $X$ and a Hermitian vector bundle $E$ over $X$ carrying a Hermitian connection $\nabla^E$, one can define  the (twisted by $E$) Dirac operator $D^E:\Gamma(S(TX)\otimes E)\rightarrow \Gamma(S(TX)\otimes E)$ to be the induced Dirac operator $D_{\partial M, +}^{\pi^*E}|_X$, in the sense of the previous subsection,  with $M=X\times [0,1]$ and $\pi:M\rightarrow X=X\times\{0\}$ being  the canonical projection.  The (reduced) $\eta$-invariant of $D^E$ is thus well-defined, without assuming that $X$ is the boundary of another spin manifold. 

The basic properties of $\overline{\eta}(D^E)$ have  indeed been studied in \cite{P7}. We mention in particular that in \cite[I, (4.2)]{P7}, it is proved by using heat kernel method that  the $\eta$-function $\eta(D^E, s)$ is holomorphic for ${\rm Re}(s)>-\frac{1}{2}$,\footnote{This result is later improved in \cite[I]{BF} to  ${\rm Re}(s)>- {2}$.} while in \cite[II, Section 4]{P7}, the intrinsic relation between the $\eta$-invariant and  the  Chern-Simons  invariant  \cite{CS} was exploited.  Moreover, by applying Theorem \ref{t2.1} to $M=X\times [0,1]$, one can prove an anomaly formula for $\overline{\eta}(D^E)$ with respect to the variations of $g^{TX}, \ g^E$ and $\nabla^E$. 

One of the most interesting outcome from the above mentioned anomaly formula occurs when one assumes that the unitary connection $\nabla^E$ is flat, that is, the curvature $R^E=(\nabla^E)^2$ vanishes. In this situation, one finds that the  following geometric quantity,
\begin{align}\label{2.13}
 \rho_E= \overline{\eta}\left(D^E\right)-{{\rm rk}(E)}\, \overline{\eta}\left(D^{{\bf C}}\right)\ \ \ {\rm mod}\ {\bf Z}
\end{align}
does not depend on $g^{TM}$. Thus we get a topological invariant for unitary representations of the fundamental group of $X$.  The remarkable {\it Index theorem for flat bundles} proved in \cite[III, (5.3)]{P7} gives a topological interpretation of this analytically constructed invariant.
This reminds us the relation between two other invariants for flat vector bundles: the Reidemeister torsion and the Ray-Singer analytic torsion, which will be discussed in the next section. 

We now turn to a beautiful observation   due to Witten \cite{W85}. Based on the reasonings in physics, Witten needed to compute the (reduced) $\eta$-invariant for Dirac operators on manifolds   fibering  over a circle. What Witten claimed in \cite{W85} is that when taking the adiabatic limit, the $\eta$-invariant for a fibered manifold over a circle is closely related to the holonomy of the Quillen determinant line bundle \cite{Q2}   over the circle,  constructed through the fiberwise Dirac operators.  This conjecture of Witten was proved independently by Bismut-Freed \cite{BF} and Cheeger \cite{Ch87}. Both proofs were reported at the ICM-1986 in Berkeley (\cite{B86a}, \cite{Ch86}). For an approach in the spirit of both \cite{BF} and \cite{Ch87},   see  \cite{DF}.

With the encouragement of the solution to the Witten conjecture, it is natural to study the $\eta$-invariant on general fibered manifolds, and this was first accomplished in the celebrated joint work of Bismut and Cheeger \cite{BC89}. 

To be more precise, let $Z\rightarrow X\stackrel{\pi}{\rightarrow}B$ be a smooth fibration of closed oriented spin manifolds, with compatible orientations and spin structures. We assume $\dim X$ is odd.   Let $T^VX$ (usually denoted by $TZ$) be the vertical tangent bundle of the fibration and $T^HX\subset TX $ be a fixed  horizontal subbundle of $TX$. Then we have the splitting $TX=T^HX\oplus TZ$. 

Let   $g^{TZ}$ be a Euclidean metric on $TZ$. Let $g^{TB}$ be a Riemannian metric on $TB$. It determines a Euclidean metric $ \pi^*g^{TB}$ on $T^HX$. Let $g^{TX}$ be the Riemannian metric on $TX$ given by the orthogonal direct sum $g^{TX}=\pi^*g^{TB}\oplus g^{TZ}$. 

For any $\varepsilon>0$, let $g_\varepsilon^{TX}$ be the rescaled metric on $TX$ defined by
\begin{align}\label{2.14}
g_\varepsilon^{TM}=\frac{\pi^*g^{TB}}{\varepsilon}\oplus g^{TZ}.
\end{align}

Let $E$ be a Hermitian vector bundle over $X$ carrying a Hermitian connection $\nabla^E$. Then for any $b\in B$, on the  fiber $Z_b=\pi^{-1}(b)$, the restrictions of $g^{TZ}$, $g^E$ and $\nabla^E$ on $Z_b$ determine a Dirac operator $D^{E}_{Z_b}$ on $Z_b$. 

Let $D_\varepsilon^E$ be the Dirac operator on $X$ associated to $g_\varepsilon^{TX}$, $g^E$ and $\nabla^E$. The following result computes the adiabatic limit of $\overline{\eta}(D_\varepsilon^E)$ when $\varepsilon$ tends to zero. 

\begin{thm}\label{t2.5}\mbox{\rm (Bismut-Cheeger \cite{BC89})} If for any $b\in B$, $D_{Z_b}^E$ is invertible, then there is $\varepsilon_0>0$ such that for any $0<\varepsilon<\varepsilon_0$,   $D_\varepsilon^E$ is invertible. Moreover, the following identity holds:
\begin{align}\label{2.15}
\lim_{\varepsilon\rightarrow 0}\overline{\eta}\left(D_\varepsilon^E\right)=\int_B\widehat{A}\left(TB,\nabla^{TB}\right)\widehat{\eta},
\end{align}
where $\nabla^{TB}$ is the Levi-Civita connection of $g^{TB}$ and $\widehat{\eta}\in\Omega^*(B)$ is the $\widehat{\eta}$-form canonically constructed in \cite{BC89}. 
\end{thm}

In the special case where $B$ is a circle, there is no specific assumptions on the fiberwise Dirac operators in the proofs \cite{BF}, \cite{Ch87} of the Witten holonomy conjecture. For general fibrations, a refinement of Theorem \ref{t2.5} was obtained by Dai \cite{D} under the assumption that $\ker(D_b^E)$, $b\in B$, form a vector bundle (denoted by $\ker(D ^E_Z)$) over $B$.  If this assumption holds, then $\ker (D ^E_Z)$ admits a naturally induced Hermitian metric as well as a Hermitian connection. Thus, one can construct the associated (twisted by $\ker (D^E_Z)$) Dirac operator $D^{\ker (D^E_Z)}$ on $B$. The main result in \cite{D} shows that  the $\widehat{\eta}$-form of Bismut-Cheeger appeared in (\ref{2.15}) is still well-defined under this condition. Moreover, the following identity holds (compare with \cite{L94}):  
\begin{align}\label{2.16}
\lim_{\varepsilon\rightarrow 0}\overline{\eta}\left(D_\varepsilon^E\right)\equiv
\overline{\eta}\left(D ^{\ker (D^E_Z)}\right)+\int_B\widehat{A}\left(TB,\nabla^{TB}\right)\widehat{\eta} \ \ {\rm mod}\ \ {\bf Z}. 
\end{align}
In particular, if $\dim B$ is even, then one has
\begin{align}\label{2.17}
\lim_{\varepsilon\rightarrow 0}\overline{\eta}\left(D_\varepsilon^E\right)\equiv
\frac{\dim \left(\ker D ^{\ker (D^E_Z)}\right)}{2}+\int_B\widehat{A}\left(TB,\nabla^{TB}\right)\widehat{\eta} \ \ {\rm mod}\ \ {\bf Z}. 
\end{align}

In \cite{BC92}, by applying   Theorem \ref{t2.5} to  flat torus bundles, a new proof of the Hirzebruch conjecture, which was first proved by Atiyah-Donnelly-Singer \cite{ADS} and M\"uller \cite{Mu84}, is given.  The $\widehat{\eta}$-form for circle bundles was computed in \cite{DZ95} and \cite{Z94}. It was applied in \cite{Z94} to obtain a higher dimensional Rokhlin type congruence which at the same time extends the divisibility result of Atiyah-Hirzebruch \cite{AH}, which states that the $\widehat{A}$-genus of an $8k+4$ dimensional closed spin manifold is an even integer, to the case of  spin$^c$-manifolds. 

When $B$ is a point, then $\widehat{\eta}$ is simply   (half of)   the  $\eta$-invariant of $D^E$. Thus, the $\widehat{\eta}$-form of Bismut-Cheeger generalizes the   $\eta$-invariant of Atiyah-Patodi-Singer to the case of families. The natural question of whether there is a families generalization of the Atiyah-Patodi-Singer index theorem \ref{t2.1} was answered positively in \cite[II]{BC90}, \cite{BC91} and \cite{MeP}. On the other hand, the concept of spectral flow was generalized to the families index theory  in  \cite{DZ98}. 

The ideas and methods in \cite{BC89}, \cite{BC90}, \cite{BC91}, \cite{MeP} and \cite{DZ98} have been further generalized to the framework of noncommutative index theory. We refer to \cite{LP}, \cite{LP04} for an overview.

The $\eta$-invariant and the $\widehat{\eta}$-form have also played important roles in the actively studied {\it differential $K$-theory}. Besides (\ref{2.15})-(\ref{2.17}), the Riemann-Roch property of the $\eta$-invariant under embedding, which is established in \cite{BZ93},  is also used in the proof of the index theorem in differential $K$-theory \cite{FL}. On the other hand,  as a simple application of the main result in \cite{BZ93}, a geometric formula for the mod ${\bf Z}$ component of the (reduced) $\eta$-invariant for Dirac operators is given in \cite{Z05}.  This later  formula, which can be proved independently, is in turn used in \cite{FXZ} to give a more geometrical proof of the embedding formula in \cite{BZ93}.

We now describe another line of applications of the $\eta$-invariant.  It concerns the heat kernel proof of the index theorem for Toeplitz operators on {\it odd} dimensional manifolds. 

Recall that $X$ is an odd dimensional closed oriented spin Riemannian manifold, $E$ a Hermitian vector bundle over $X$ carrying a Hermitian connection, and $D^E:\Gamma(S(TX)\otimes E)\rightarrow \Gamma(S(TX)\otimes E)$ is the Dirac operator which is elliptic and formally self-adjoint. 

Now let $N>0$ be an integer and ${\bf C}^N_X$ be a trivial vector bundle over $X$ carrying the  canonical trivial Hermitian metric and connection. Then $D^E$ extends naturally to a Dirac operator acting on $\Gamma(S(TX)\otimes E\otimes {\bf C}^N_X)$, which we still denote $D^E$. 

Let $g\in\Gamma({\rm Aut}({\bf C}^N_X))$ be a fiberwise automorphism of ${\bf C}^N_X$. It extends to a fiberwise automorphism ${\rm Id}_{S(TX)\otimes E}\otimes g$, which we still denote by $g$. Then $g$ acts as a bounded linear operator on the $L^2$-completion $L^2(S(TX)\otimes E\otimes {\bf C}^N_X)$ of $\Gamma(S(TX)\otimes E\otimes {\bf C}^N_X)$.
Let $P^E_{\geq 0}:L^2(S(TX)\otimes E\otimes {\bf C}^N_X)\rightarrow L^2_{\geq 0}(S(TX)\otimes E\otimes {\bf C}^N_X)$ be the orthogonal projection defined as in between  (\ref{2.4}) and (\ref{2.5}). 

We define  the {\it Toeplitz operator} $T_g^E$ by
\begin{align}\label{2.18}
T_g^E=P^E_{\geq 0}gP^E_{\geq 0}: L^2_{\geq 0}\left(S(TX)\otimes E\otimes {\bf C}^N_X\right)\longrightarrow L^2_{\geq 0}\left(S(TX)\otimes E\otimes {\bf C}^N_X\right).
\end{align}

It was observed in \cite{BD} that $T_g^E$ is a zeroth order elliptic pseudodifferential operator, so that one can apply the Atiyah-Singer index theorem \cite{ASI} to compute its index. More precisely, one gets
\begin{align}\label{2.19}
{\rm ind}\left( T_g^E\right)=-\left\langle\widehat{A}(TX){\rm ch}(E){\rm ch}(g),[X]\right\rangle,
\end{align}
where ${\rm ch}(g)$ is the odd Chern character of $g$, when viewing $g$ as an element in $K^1(X)$ (cf. \cite[Chap. 1]{Z01}). 

Following \cite{BoW93}, we   outline below a proof of (\ref{2.19}) using heat kernels and the $\eta$-invariant. The first observation is that by the homotopy invariance of the index, one may deform $g$ to assume that it is fiberwise unitary. Now as $g$ is fiberwise unitary, we see that the operator $g^{-1}D^Eg:\Gamma(S(TX)\otimes E\otimes {\bf C}^N_X)\rightarrow \Gamma(S(TX)\otimes E\otimes {\bf C}^N_X)$ is formally self-adjoint. Moreover,  $D_t^E=(1-t)D^E+tg^{-1}D^Eg$,   $t\in [0,1]$, form a smooth family of formally self-adjoint elliptic operators acting on the same space $\Gamma(S(TX)\otimes E\otimes {\bf C}^N_X)$.  

One first observes that the   index of   $T^E_g$ can be identified with the spectral flow for the above family $\{D_t^E:0\leq t\leq 1\}$. That is, one has
\begin{align}\label{2.20}
{\rm ind}\left( T_g^E\right)=-{\rm sf}\left\{D_t^E:0\leq t\leq 1\right\}. 
\end{align}
On the other hand,  by the variation formula for the $\eta$-invariant as in \cite[III]{P7}, one has
\begin{align}\label{2.21}
\overline{\eta}\left(D_{t=1}^E\right)-\overline{\eta}\left(D_{t=0}^E\right)={\rm sf}\left\{D_t^E:0\leq t\leq 1\right\}+\int_0^1\frac{d\overline{\eta}\left(D_t^E\right)}{dt}\,dt.
\end{align} 
Now it is clear that $\overline{\eta} (D_{t=1}^E )=\overline{\eta} (g^{-1}D ^E g)=\overline{\eta} ( D ^E )=\overline{\eta}\left(D_{t=0}^E\right)$, from (\ref{2.20}) and (\ref{2.21}) one gets
\begin{align}\label{2.22}
{\rm ind}\left( T_g^E\right)=\int_0^1 \frac{d\overline{\eta}\left(D_t^E\right)}{dt}\,dt,
\end{align} 
with the right hand side being able to be evaluated by local index computations (i.e., the heat kernel method), which leads to (\ref{2.19}).    

Inspired by Theorem \ref{t2.1}, it is natural to ask whether (\ref{2.19}) can be generalized to the case of manifolds with boundary. This problem was first studied in \cite{DoW}, where an index formula is obtained, under the assumption that the restriction of $g$ on the boundary is an identity. No boundary contribution appears in the index formula under this assumption. 

On the technical side, \cite{DoW}, which is inspired by an earlier note of Singer \cite{S87}, initiated  the study of $\eta$-invariants on manifolds with boundary, which have been studied actively since then. Here we only mention the papers \cite{Bu}, \cite{KL} and \cite{Mu94} from the vast literature on this subject. 

A solution to the problem of   generalizing   (\ref{2.19}) to the case of manifolds with boundary, without any condition on $g$, is  later given  in \cite{DZ06}, where the results in \cite{KL} and \cite{Mu94} play important roles. The result proved in \cite{DZ06} may be thought of as an odd dimensional analogue of Theorem \ref{t2.1}. In particular, there appears a boundary contribution in the index formula established in \cite{DZ06}, which is indeed an $\eta$-type invariant.

Certainly there are many other implications  of Theorem \ref{t2.1} and the $\eta$-invariant which are not touched in this very brief account. We refer to the  survey papers  \cite{Mu94a} and \cite{Go}  for some of the topics we did not discuss  above.

\section{Reidemeister torsion and Ray-Singer analytic torsion} \label{s2a}

As we have pointed out in Introduction, Patodi's work touches  almost all aspects of the index theory. Besides the works we have described in the previous two sections, in his ICM talk \cite{P8}, his interests moved on to two other important problems: the first is  to give an explicit combinatorial formula for the Pontryagin classes; the second is to give an analytic interpretation of the Reidemeister torsion  \cite{R35}. 

The former problem is still open up to now, while Patodi   made his own contributions in \cite{P9}.  Here we only mention  that the famous Chern-Simons invariant \cite{CS} was in fact originated from an attempt to give a purely combinatorial expression for the first Pontryagin class of 4-manifolds. 

In this section, we will concentrate on the other problem touched in \cite{P8}, which  concerns the Reidemeister torsion.

The Reidemeister torsion  is a classical invariant in topology associated
to   orthogonal  representations of the fundamental group of a CW
complex.

Inspired by the Atiyah-Singer index theory, Ray and Singer \cite{RS71} studied
  an analogue of the Reidemeister torsion for the de Rham
complex on a smooth manifold. They called this analogue the
analytic torsion and discovered that it has a lot of  
properties similar to  those of the Reidemeister torsion. They further made
the conjecture that their analytic torsion (now widely called the
Ray-Singer torsion) equals to the Reidemeister torsion.

For simplicity, we start directly with a smooth closed manifold $M$ which we assume to be oriented. Let $\rho:\pi_1(M)\rightarrow O(N)$ be an orthogonal  representation of the fundamental group of $M$. Then $\rho$ determines an orthogonal  flat vector bundle $E_\rho=\widetilde{M}\times_\rho{\bf R}^N$ over $M$, where $\widetilde{M}$ is the universal covering of $M$. 

Take any triangulation   $\{c_\alpha\} $ of $M$, where each $c_\alpha$ is a simplex. Over each  simplex $c_\alpha$, $E_\rho$ can be trivialized tautologically. Thus the boundary operator $\partial: c_\alpha\rightarrow \partial c_\alpha$ extends naturally to a twisted boundary operator $\partial_\rho: c_\alpha\otimes ( E_\rho|_{c_\alpha}) \rightarrow (\partial c_\alpha)\otimes( E_\rho|_{\partial c_\alpha})$.

For any simplex $c_\alpha$, let $e_{\alpha 1},\,\cdots,\,e_{\alpha N}$ be the trivialized global orthonormal basis of $E_\rho|_{c_\alpha}$, let $[c_\alpha\otimes e_{\alpha i}]$, $1\leq i\leq N$, denote the real line generated by $c_\alpha\otimes e_{\alpha i}$. For any integer $0\leq q\leq \dim M$, let $C_q$ denote the vector space obtained by the direct sum $C_q=\oplus_{i=1}^N\oplus_{ \dim c_\alpha=q}[c_\alpha\otimes e_{\alpha i}]$. With these data one can  construct a combinatorial complex $(C_*,\partial_\rho)$ with coefficient $E_\rho$,
\begin{align}\label{3.1}
\left(C_*,\partial_\rho\right): 0\longrightarrow C_{\dim M}\stackrel{\partial_\rho}{\longrightarrow} C_{\dim M-1}
\stackrel{\partial_\rho}{\longrightarrow}\cdots\stackrel{\partial_\rho}{\longrightarrow}C_0\longrightarrow 0,
\end{align} 
i.e., one has $\partial_\rho^2=0$. 
For simplicity, we also assume that the homology of the complex $(C_*,\partial_\rho)$ vanishes:  $H_*(C_*,\partial_\rho)=\{0\}$. Then for any $0\leq q\leq \dim M$, $\ker (\partial_\rho|_{C_q})\simeq {\rm Im} (\partial_\rho|_{C_{q+1}})$.\footnote{We set $C_{-1}=C_{\dim M+1}=\{0\}$. }

For any $0\leq q\leq \dim M$, let $b_{qj}$'s be a basis of ${\rm Im} (\partial_\rho|_{C_{q+1}})$, let $b_{(q-1)j}'\in C_{q}$  be  such that $\partial_\rho b_{(q-1)j}'=b_{(q-1)j}$.  Then the $(b_{qj}, b_{(q-1)j}')$'s form a basis of $C_q$. Recall that the $c_\alpha\otimes e_{\alpha i}$'s with $\dim c_\alpha=q$, $1\leq i\leq N$, form another basis $c_q$ of $C_q$. Let $[b_{qj}, b_{(q-1)j}' : c_q]>0$ be the absolute value of the determinant with respect to these two basis. Then one sees that this determinant depends only on the $b_{(q-1)j}$'s, not on the lifts  $b_{(q-1)j}'$'s. Thus we may denote it by $[b_{qj}, b_{(q-1)j}: c_q]$. 

One can now define the {\it Reidemeister torsion} $\tau(C_*,\partial_\rho)$ of the complex $(C_*,\partial_\rho)$ by the following formula, 
\begin{align}\label{3.2}
\log \tau\left(C_*,\partial_\rho\right)=\sum_{q=0}^{\dim M} (-1)^q \log \left[b_{qj}, b_{(q-1)j}: c_q\right]. 
\end{align} 
It is indeed a well-defined quantity as one can show that the right hand side of (\ref{3.2}) does not depend on the choices of the $b_{qj}$'s.   One can show further that $\tau(C_*,\partial_\rho)$ is invariant under subdivisions of the triangulation. Thus, it does not depend on the triangulations in its definition and is a well-defined combinatorial invariant. We denote it by $\tau_\rho$ and call it the Reidemeister torsion associated to $\rho$. 

Ray and Singer observed that the Reidemeister torsion defined in (\ref{3.2}) admits an expression through combinatorial Laplacians. To be more precise, let $C_*$ carry a Hermitian metric such that the $ c_\alpha\otimes e_{\alpha i}$'s with $0\leq\dim c_\alpha\leq \dim M$, $1\leq i\leq N$, form an orthonormal basis. Let $\partial_\rho^*:C_*\rightarrow C_*$ be the adjoint of $\partial_\rho$ with respect to this Hermitian metric, then for any $0\leq q\leq \dim M$, $\partial_\rho^*:C_q\rightarrow C_{q+1}$. 

Let $\Box^{C_*} _\rho:C_*\rightarrow C_*$ be the Laplacian defined by $\Box^{C_*}_\rho =\partial_\rho^*\partial_\rho+\partial_\rho\partial_\rho^*$. For any $0\leq q\leq \dim M$, let $\Box^{C_*}_{  q,\rho}:C_q\rightarrow C_q$ be the restriction of $\Box^{C_*}_\rho $ on $C_q$. Since we have assumed that $H_*(C_*,\partial_\rho)=\{0\}$, we know that each $\Box^{C_*}_{  q,\rho}$ is invertible.  The following identity is proved in \cite[Proposition 1.7]{RS71}, 
\begin{align}\label{3.3}
\log \tau_\rho=\frac{1}{2}\sum_{q=0}^{\dim M} (-1)^{q+1}q \log\left(\det \left(\Box^{C_*}_{ q,\rho}\right)\right).
\end{align} 

Inspired by (\ref{3.3}), Ray and Singer defined in \cite{RS71} an analogue for the de Rham complex with coefficient  $E_\rho$, i.e., the following complex where for any $0\leq q\leq \dim M$,  $\Omega^q(M,E_\rho)=\Gamma(\Lambda^q(T^*M)\otimes E_\rho)$,
\begin{align}\label{3.4}
\left(\Omega^*(M,\rho),d_\rho\right): 0\rightarrow\Omega^0(M,E_\rho)\stackrel{d_\rho}{\rightarrow}\Omega^1(M,E_\rho)\stackrel{d_\rho}{\rightarrow}\cdots\stackrel{d_\rho}{\rightarrow}\Omega^{\dim M}(M,E_\rho)\rightarrow 0, 
\end{align}
where $d_\rho: \Omega^*(M,E_\rho)\rightarrow  \Omega^*(M,E_\rho)$ is the (twisted) exterior differential operator with coefficient $E_\rho$. 

Let $g^{TM}$ be a Riemannian metric on $TM$. 

By using $g^{TM}$ and the orthogonal  flat structure on $E_\rho$, one can produce a canonically induced inner product   on $ \Omega^*(M,E_\rho)$. Let $d^*_\rho: \Omega^*(M,E_\rho)\rightarrow  \Omega^*(M,E_\rho)$ be the formal adjoint of $d_\rho$ with respect to this inner product. Let $\Box_\rho=d_\rho^*d_\rho+d_\rho d_\rho^*$ be the associated Laplacian. For any $0\leq q\leq \dim M$, let $\Box_{q,\rho}:\Omega^q(M,E_\rho)\rightarrow \Omega^q(M,E_\rho)$ be the restriction of $\Box_\rho$ on $\Omega^q(M,E_\rho)$. 

If one wants to define an analogue of (\ref{3.3}), one needs to at first define a kind of ``determinant'' for the Laplacians $\Box_{q,\rho}$'s. The obvious difficulty is that the $\Box_{q,\rho}$'s act on infinite dimensional spaces, so one can not define a  ``determinant'' directly for them. The beautiful observation in \cite{RS71} is that one can use a $\zeta$-function regularization to define  a kind of {\it regularized determinant} which in the finite dimensional case coincides with the usual determinant. 

Since we have assumed $H_*(C_*,\partial_\rho)=\{0\}$, by using the Hodge theorem one knows that for each $0\leq q\leq \dim M$, $\Box_{q,\rho}$ is invertible. Following \cite{RS71}, by standard elliptic theory, for any $s\in {\bf C}$ with ${\rm Re}(s)>>0$, the $\zeta$-function in the following formula is well-defined:
\begin{align}\label{3.5}
 \zeta_{q,\rho}(s)=\sum_{\lambda\in{\rm Spec}(\Box_{q,\rho})}\frac{1}{\lambda^s} .
\end{align}
It is easy to show that $\zeta_{q,\rho}(s)$ can be extended to a meromorphic function on ${\bf C}$ which is holomorphic at $s=0$. 

Following \cite[Definition 1.6]{RS71}, we now define the {\it Ray-Singer analytic torsion} $T_\rho$ by the following formula,
\begin{align}\label{3.6}
 \log T_\rho=\frac{1}{2}\sum_{q=0}^{\dim M}(-1)^qq \left.\frac{d \zeta_{q,\rho}(s) }{ds}\right|_{s=0}. 
\end{align} 
By using (\ref{3.3}), one sees easily that $\tau_\rho$ admits a similar formula in terms of the combinatorial Laplacians. To be more precise, if we set $\zeta(\Box^{C_*}_{q,\rho},s)=\sum_{\lambda\in{\rm Spect}(\Box^{C_*}_{q,\rho})}\frac{1}{\lambda^s}$ for $0\leq q\leq\dim M$, then we can rewrite (\ref{3.3}) as 
\begin{align}\label{3.7}
  \log \tau_\rho=\frac{1}{2}\sum_{q=0}^{\dim M}(-1)^qq \left.\frac{d \zeta\left(\Box^{C_*}_{q,\rho},s\right) }{ds}\right|_{s=0}. 
\end{align}

It is proved in \cite{RS71} that when $\dim M$ is even, $T_\rho=1$, while when $\dim M$ is odd, $T_\rho$ does not depend on $g^{TM}$, i.e., it is a smooth invariant! Ray and Singer also proved in \cite{RS71} that many properties of $\tau_\rho$ are also satisfied by $T_\rho$. They raised the conjecture that these two invariants actually equal, i.e., the following identity holds,\footnote{The original conjecture also covers   the case where $H_*(C_*,\partial_\rho)\neq\{0\}$.}
\begin{align}\label{3.8}
   T_\rho=\tau_\rho. 
\end{align} 

It is this conjecture which Patodi dealt with in his ICM talk \cite{P8}, and also in his joint paper with Dodziuk \cite{P10}. In fact, as we have mentioned, $\tau_\rho$ does not depend on subdivisions of the triangulation in its definition. So one of the ways one could expect to prove (\ref{3.8}), in view of 
(\ref{3.6}) and (\ref{3.7}), is that if one keeps doing the subdivision, then the $\zeta$-function associated to the combinatorial Laplacians might converge in some sense to the $\zeta$-function of the Laplacians of the de Rham-Hodge operator. Moreover, if this would hold,  it might lead to a proof of (\ref{3.8}). 

The study of the approximation of the combinatorial Laplacian under subdivisions were studied by Dodziuk in his thesis \cite{D76}.  The method in \cite{D76} was then further developed in the joint paper \cite{P10},  of which some of the main results were announced in \cite{P8}, with the purpose of proving the Ray-Singer conjecture. While \cite{P10} does not give a proof of (\ref{3.8}), the methods and results obtained are inspiring. We just mention one approximation result as follows. 

In \cite[Theorem 4.10]{P10}, it is stated that for any compact subset $U$ in the half plane ${\rm Re}(s)>\frac{\dim M}{2}$, there is a sequence of triangulations $c_{\alpha, k}$, such that for any $0\leq q\leq \dim M$, the corresponding sequence of $\zeta$-functions  $\zeta (\Box^{C_{*,k}}_{q,\rho},s )$ converges to $\zeta_{q,\rho}(s)$ uniformly on $U$.  

While the above result does not lead to (\ref{3.8}) which involves informations of the $\zeta$-functions around $s=0$, it is obviously  an important  step towards a better understanding of (\ref{3.8}). In fact, a program was proposed in \cite{P10} towards  a possible  proof of (\ref{3.8}) by stating two further conjectures concerning the approximations of combinatorial Laplacians.

The   Ray-Singer conjecture (\ref{3.8})  was finally proved   independently in the celebrated papers of 
  Cheeger \cite{Ch79}  and M\"uller \cite{Mu}, by different methods. While Cheeger's method is based mainly on surgery, what M\"uller did is indeed to further develop the combinatorial approximation methods in \cite{P10}, by completing the program proposed in \cite{P10}. Thus, we could say that in some sense Patodi  also played a pioneering role in the study  of the Ray-Singer conjecture, now the Cheeger-M\"uller theorem. 

As it is stated, the Cheeger-M\"uller theorem holds for
orthogonal  representations of the fundamental group of smooth manifolds.   The natural question of whether (\ref{3.8}) can be extended to other representations of the fundamental group was solved by 
M\"uller \cite{Mu93} for unimodular representations and by Bismut-Zhang   \cite{BZ92} for arbitrary representations. 
 
The method used in \cite{Mu93} extends those of Cheeger's in \cite{Ch79}. The method developed in \cite{BZ92}
 is purely
analytic and makes use of  the Witten deformation \cite{W82} for the de Rham complex
by a Morse function. Moreover, the local index theoretic technique, as developed from \cite{P1}, plays essential roles in \cite{BZ92}. 

Ray-Singer also defined analytic torsion for manifolds with boundary in \cite{RS71}. This concept has played important roles in the proofs of Cheeger and M\"uller on the Ray-Singer conjecture (\ref{3.8}). The boundary conditions appear here are of local type, in comparing with the global boundary conditions discussed in Section \ref{s2}.  The Ray-Singer torsion for manifolds with boundary has   been studied extensively. Here we only mention two  papers  of Br\"uning-Ma (\cite{BrM06} and \cite{BrM13}),  which deal with  the most general case of arbitrary flat vector bundles as well as arbitrary metrics near the boundary. In particular, a local Gauss-Bonnet-Chern theorem, generalizing the work of Patodi \cite{P1} to the case of manifolds with boundary, is established in \cite{BrM06}, where no condition on   metrics near the boundary is assumed. 

Inspired by the construction of the Bismut-Cheeger $\widehat{\eta}$-form, which we have described in the last section,  Bismut-Lott \cite{BL95} constructed for smooth fibrations (real) analytic torsion forms   of which the degree zero term is   the fiberwise Ray-Singer torsion. A vast generalization of the Cheeger-M\"uller and Bismut-Zhang theorems, to the Bismut-Lott torsion form for fibrations, was proved by Bismut-Goette in \cite{BG}. For further developments concerning the torsion forms, we refer to the survey paper \cite{Go09}. 

Recall that in (\ref{2.13}), we have described the Atiyah-Patodi-Singer $\rho$-invariant for flat vector bundles through  $\eta$-invariants. It is natural to ask whether this invariant is related to the Ray-Singer analytic torsion, which is also a topological invariant for flat vector bundles. In the Introduction to    \cite[I]{P7}, it is mentioned that a unification  of the $\eta$-invariant and the Ray-Singer torsion was not known. While this is still the case up to now, the fact that they both appear prominently  in the Chern-Simons gauge theory   (cf. \cite{W89} and \cite{BW91}), which reflects other deep implications of the $\eta$-invariant and the analytic torsion,  may shed a light on this mystery (cf. \cite{MZ06} where a holomorphic function on the representation space of fundamental group, with absolute value being equal to the associated Ray-Singer analytic torsion, is constructed by using generalized versions of the above $\eta$-invariant and analytic torsion). 
On the other hand, relations between the $\widehat{\eta}$-form and the Bismut-Lott torsion form have been exploited in \cite{MZ8} and \cite{DZ10}.

In another direction, Ray and Singer also defined analytic torsion for the Dolbeault complex on complex manifolds in  \cite{RS73}, which we now call the Ray-Singer holomorphic torsion. While this holomorphic torsion is no longer a topological invariant, it has played   important roles in complex geometry,  arithmetic geometry and mathematical physics. We refer to Bismut's ICM-1998 Plenary Lecture \cite{B98} for a survey on some of the developments concerning the holomorphic torsion.

As for the most significant recent advance of the (local) index theory, we mention the groundbreaking theory of hypoelliptic Laplacians developed by Bismut and his collaborators  in \cite{B05}, \cite{B08a}, \cite{BL08},  \cite{B11} and \cite{B13}, of which two surveys are given in \cite{B08} and \cite{B12}.  The motivation of this theory comes in part from an attempt to generalize the methods and results developed in \cite{BZ92} to the loop space of the underlying smooth manifold.\footnote{See also \cite{S14} for a proof  via   hypoelliptic Laplacian of the main result in \cite{BZ92}.} In this whole new development, we find that the local index techniques, the  Atiyah-Patodi-Singer $\eta$-invariant and the Ray-Singer analytic torsion all find their new places. We believe that if Patodi would still  be alive, he would be happy to see that the local index theory he helped to create   has become one of the central  areas in global analysis.

  \def\cprime{$'$} \def\cprime{$'$}
\providecommand{\bysame}{\leavevmode\hbox to3em{\hrulefill}\thinspace}
\providecommand{\MR}{\relax\ifhmode\unskip\space\fi MR }
\providecommand{\MRhref}[2]{%
 \href{http://www.ams.org/mathscinet-getitem?mr=#1}{#2}
}
\providecommand{\href}[2]{#2}

$\ $

\noindent Chern Institute of Mathematics \& LPMC 

\noindent Nankai University

\noindent Tianjin 300071

\noindent P. R. China

$\ $

\noindent {\it Email: weiping@nankai.edu.cn}

\end{document}